\documentstyle[11pt,leqno]{article}
\def\secsep{\mbox{}\\[1\baselineskip]\mbox{}}

\def\sect#1{\secsep{\large\bf #1}}

\begin{document}
\parindent    = 25pt
\baselineskip = 16pt

\mbox{}\\
\begin{center}
{\Large\bf The number of ramified covering of a Riemann surface\\ by Riemann 
surface}
\vskip 0.4in

An-Min Li\footnote{Partially supported by a NSFC and a QiuShi 
grant. \quad e-mail address: amli@scu.edu.cn }\qquad
Guosong Zhao\footnote{Partially supported by a NSFC.
\quad e-mail address: gszhao@scu.edu.cn}\qquad
Quan Zheng\footnote{Partially supported by a postdoctor foundation of 
Dalian
University of Technology and a Youth Foundation of Sichuan 
University.}\qquad

(Department of Mathematics, Sichuan University\\
610064, Chengdu, Sichuan, PRC)

(July 27, 1999)\\
\vskip 0.4in 
{\bf Abstract}\\[\baselineskip]

\parbox{14cm}{\quad
 Interpreting the number of ramified covering of a Riemann surface by Riemann 
surfaces as the  relative Gromov-Witten invariants and
applying a gluing formula, we derive a recursive
formula for the number of ramified covering of a Riemann surface by Riemann surface 
with elementary branch points and prescribed ramification type 
over a special point.
}
\end{center}

\sect{1\quad Introduction}
\setcounter{section}{1}  \\

Let $\Sigma^g$ be a compact connected Riemann surface of genus $g$ and $\Sigma^h$
 a compact connected Riemann surface of genus $h$ $(g \geq h\geq 0).$ 
A ramified
covering of $\Sigma^h$ of degree $k$ by $\Sigma^g$ is a non-constant 
holomorphic map
$f:\Sigma^g\to \Sigma^h$ such that $|f^{-1}(q)|=k$ for all but a finite 
number of points
$q\in \Sigma^h$, which are called branch points. Two ramified coverings $f_1$ 
and $f_2$
 are said to be equivalent if there is a 
homeomorphism
$\pi: \Sigma^g\to \Sigma^g$ such that $f_1=f_2\circ\pi$. A ramified covering 
$f$ is called almost
simple if $|f^{-1}(q)|=k-1$ for each branch point but one, that is 
denoted by $\infty$.
 If 
$\alpha_1,\ldots,\alpha_m$ are
the orders of the preimage of $\infty$, then the ordered $m$-tiple pair
$(\alpha_1,\cdots,\alpha_m)=\alpha$ is a partition of $k$, denoted
by $\alpha \vdash k,$
 and is called the ramification type of $f$ (at $\infty$). We call 
$m$ the
length of $\alpha$, denoted by $l(\alpha)=m$.
Let $\mu_{h,m}^{g,k}(\alpha)$ be the number of equivalent almost simple
 covering of $\Sigma^h$ by $\Sigma^g$ with ramification type of $\alpha$.
 How to determine $\mu_{h,m}^{g,k}(\alpha)$ 
is known as the
Hurwitz Enumeration Problem. It is Hurwitz who first gave a
explicit expression for $\mu_{0,m}^{0,k}(\alpha)$, see [H]. $\mu_{h,m}^{g,k}(\alpha)$ is
called Hurwitz number.
Many  mathematicians contribute to this problem for the case $h = 0$. J.D\'{e}nes [D]
gave a formula for $g=0,l(\alpha)=1$, and V.I.Arnol'd [A] for $g=0,
l(\alpha)=2$. By a combinatorial method, I.P. Goulden and D.M. Jackson gave an explicit
formula for $g= 0,1,2.$ (see [GJ1,2,3,], [GJV]). Physicists M.Crescimanno and W.Taylor [CT] solved
the case when $g=0$, and $\alpha$ is the identity. R.Vakil [V] derived explicit
expressions for $g=0,1$,  using a deformation theory of algbraic geometry . \\

In this paper we interpret Hurwitz number $\mu_{h,m}^{g,k}(\alpha)$ as the  relative
Gromov-Witten invariants defined by  Li and  Ruan [LR], then applying
a gluing formula we derive a recursive formula for $\mu_{h,m}^{g,k}(\alpha)$ with $ g \geq h.$\\

Suppose $\alpha =(\alpha_1,\cdots,\alpha_m)$, put $\tilde{J}(\alpha)=\left\{
(\alpha_i,\alpha_j) 
| 1\leq i<j\leq m\right\}/\sim$, where $(\alpha_i, \alpha_j)$ $\sim$ $(\alpha_s, \alpha_t)$
 iff $(\alpha_i,\alpha_j)=(\alpha_s,\alpha_t)$, or 
$(\alpha_i,\alpha_j)=(\alpha_t,\alpha_s)$. For every equivalent class
    $[(\alpha_i, \alpha_j)]$ of $\tilde{J}(\alpha)$
we choose a representative element, say $(\alpha_i, \alpha_j)$, and associate
it 
 a ordered $(m-1)$-tuple 
 $\theta=
(\alpha_1,\cdots,\hat{\alpha}_i,\cdots,\hat{\alpha}_j,\cdots$,
$\alpha_m$,$\alpha_i+\alpha_j)$, where the caret means that the term is
 omitted.
Then we obtain a  set $J(\alpha)$ of ordered $(m-1)$-tuples.
 For $\theta=(\alpha_1, \cdots, \hat{\alpha}_i,$ $ \cdots, \hat{\alpha}_j,
 \cdots, \alpha_m,$ $ \alpha_i+\alpha_j)$ $ \in J(\alpha)$, we define a
  integer
\begin{equation}
I_1(\theta) =\left\{
	\begin{array}{ccl}
        \frac{1}{2}(\alpha_i+\alpha_j)\cdot\#\{\lambda\in\theta|
        \lambda=\alpha_i+\alpha_j\}  & \mbox{ if } & \alpha_i=\alpha_j,\\
        (\alpha_i+\alpha_j)\cdot\#\{\lambda\in\theta|
        \lambda=\alpha_i+\alpha_j\} & \mbox{ if } & \alpha_i\not=\alpha_j.
	\end{array}
\right.\\
\end{equation}
 For every $\alpha_i$ in $\alpha$, we construct a set $ C_{\alpha_i}(\alpha)$ =
$\left\{\omega=(\alpha_1,\cdots,\hat{\alpha}_i,
\cdots,\alpha_m,\rho,\alpha_i-\rho)\ |\
 1\leq\rho\leq [\frac{\alpha_i}{2}] \right\}$ of ordered $(m+1)$-tiples.
 Let $\alpha_{1^\prime}, \alpha_{2^\prime}, \cdots, \alpha_{l^\prime}$
 be all  elements in $\alpha$ with distinct value. Put 
 $C(\alpha)$ =
 $C_{\alpha_{1^\prime}}(\alpha)\cup C_{ \alpha_{2^\prime}}(\alpha)\cup \cdots
 \cup C_{\alpha_{l^\prime}}(\alpha).$
 For every $\omega=(\alpha_1,\cdots,\hat{\alpha}_i,\cdots,\alpha_m,$
 $\rho,\alpha_i-\rho)$
$ \in C(\alpha)$, we associate it a number 
\begin{equation}
I_2(\omega) = \left\{
	\begin{array}{ccl}
        \frac{1}{2}\rho(\alpha_i-\rho)\cdot\#\{\lambda\in\omega|\lambda=\rho\}
        \cdot(\#\{\mu\in\omega|\mu=\alpha_i-\rho\}-1)
          & \mbox{ if } & \rho= \alpha_i-\rho,\\
        \rho(\alpha_i-\rho)\cdot\#\{\lambda\in\omega|\lambda=\rho\}
        \cdot\#\{\mu\in\omega|\mu=\alpha_i-\rho\}
                               & \mbox{ if } & \rho\not= \alpha_i-\rho.
 	\end{array}
\right.\\
\end{equation}

Dividing $\{1,2,\cdots,\hat{i},\cdots,m\}$ into two parts  $\pi_1,\pi_2$,
correspondingly, we divide $\omega =(\alpha_1,\cdots,\hat{\alpha_i},\cdots,
$\\$\alpha_m,\rho,\alpha_i-\rho)$ into two parts in  forms:
 $\omega_{\pi_1}  =
(\alpha_{\pi_1},\rho),$ $\omega_{\pi_2}=(\alpha_{\pi_2},\alpha_i-\rho).$ For example
if $\pi_1=\{1\},$ $\pi_2=\{2,\cdots,\hat{i},\cdots,m\}$, then $\omega_{\pi_1}=
(\alpha_1,\rho),\omega_{\pi_2}=(\alpha_2,\cdots,\hat{\alpha_i},\cdots,\alpha_m
,\alpha_i-\rho)$. Note that $\pi_1,\pi_2$ may be empty. Write ${\cal{P}}_\omega
= \{\pi = (\omega_{\pi_{1}},\omega_{\pi_{2}})\}/\sim,$
where  $\pi = (\omega_{\pi_1}, \omega_{\pi_2})$=$(\cdots, \rho,$ $ \cdots,
 \alpha_i-\rho)$$\sim$
 $\tilde{\pi} = (\tilde{\omega}_{\pi_1}, \tilde{\omega}_{\pi_2})$
= $(\cdots, \tilde{\rho},$ $ \cdots,  \alpha_i-\tilde{\rho})$
$\in{\cal P}_\omega$  iff $\rho = \tilde{\rho}$ and $\omega_{\pi_1}$
 and $\tilde{\omega}_{\pi_1} $ are same through a permutation. We also use
 $\pi=(\omega_{\pi_1}, \omega_{\pi_2})$ to denote its equivalent class when
 there is no confusion. For every equivalent class
 $\pi = (\omega_{\pi_1}, \omega_{\pi_2})$ $\in{\cal P}_\omega,$ we associate it
a number
\begin{equation}
I_3(\pi)=\left\{\begin{array}{ccl}
        \frac{1}{2}\rho(\alpha_i-\rho)\cdot\#\{\lambda\in\omega_{\pi_1}|
     \lambda=\rho\}        \cdot\#\{\mu\in\omega_{\pi_2}|\mu=\alpha_i-\rho\}
          & \mbox{ if } & \rho= \alpha_i-\rho,\\
        \rho(\alpha_i-\rho)\cdot\#\{\lambda\in\omega_{\pi_1}|\lambda=\rho\}
        \cdot\#\{\mu\in\omega_{\pi_2}|\mu=\alpha_i-\rho\}
                               & \mbox{ if } & \rho\not= \alpha_i-\rho.
 	\end{array}
\right.\\
\end{equation}      
In this paper, we will prove the following theorem, see Section 3:\\

{\bf  Theorem A}\quad
{\it
All $\mu_{h,m}^{g,k}(\alpha)$ can be determined by a recursive
formula:

\begin{eqnarray}
\mu_{h,m}^{g,k}(\alpha)
     &=& \sum_{\theta\in 
J(\alpha)}\mu_{h,m-1}^{g,k}(\theta)\cdot I_1(\theta)
     + \sum_{\omega\in 
C(\alpha)}\mu_{h,m+1}^{g-1,k}(\omega)\cdot I_2(\omega)\\
     & & + \sum_{\omega\in C(\alpha)}
           \sum_{\mbox{\scriptsize $\begin{array}{c} g_1+g_2=g \\ 
g_1,g_2\geq 0
                                          \end{array}$}} 
           \sum_{\mbox{\scriptsize $\begin{array}{c} k_1+k_2=k \\ 
k_1,k_2\geq 1
					  \end{array}$}}
           \sum_{\mbox{\scriptsize $\begin{array}{c} m_1+m_2=m+1\\ 
m_1,m_2\geq 1
					  \end{array}$}}
           \sum_{\pi\in {\cal P}_{\omega}} \nonumber\\
     & &    {k+m-2kh-3+2g\choose k_1+m_1-2k_1h-2+2g_1}
            \mu_{h,m_1}^{g_1,k_1}(\omega_{\pi_1}) \cdot
            \mu_{h,m_2}^{g_2,k_2}(\omega_{\pi_2})\cdot
             I_3(\pi). \nonumber
\end{eqnarray}
where  $m_i=l(\omega_{\pi_i}), i=1,2$.
}\\

Let $p=(p_1,p_2,p_3,\cdots,)$ be indeterminates, and 
$p_\alpha=p_{\alpha_1}
\cdots p_{\alpha_m}$ for $\alpha=(\alpha_1,\cdots,\alpha_m)$. Introduce a generating function
for $\mu_{h,m}^{g,k}(\alpha)$
\begin{equation}
\Phi_h(u,x,z,p)=
           \sum_{\mbox{\scriptsize $\begin{array}{c} k,m\geq 1\\ g\geq 0
					  \end{array}$}}
           \sum_{\mbox{\scriptsize $\begin{array}{c} \alpha\vdash k\\ 
l(\alpha)=m
					  \end{array}$}}
           \mu_{h,m}^{g,k}(\alpha)\cdot
            \frac{u^{k+m-2kh-2+2g}}{(k+m-2kh-2+2g)!}\cdot
                \frac{x^k}{k!}\cdot z^g p_\alpha.\nonumber
\end{equation}
Then after symmetrizing $\alpha_i$ and $\alpha_j$ in the variable
$\theta=(\alpha_1, \cdots, \hat{\alpha}_i,$ $ \cdots, \hat{\alpha}_j,
 \cdots, \alpha_m,$ $ \alpha_i+\alpha_j)$ $ \in J(\alpha)$, and $\rho$ and $\alpha_i-\rho$
 in the variable  $\omega$, $\gamma =(\alpha_1,\cdots,\hat{\alpha_i},\cdots,
$ $\alpha_m,\rho,\alpha_i-\rho)$ in (1.4), we have the following
\\

{\bf  Theorem B}\quad
{\it
The recursive formula (1.4) is equivalent to the the following partial differential equation :
\begin{equation}
\frac{\partial\Phi_h }{\partial u} 
        = \frac{1}{2}\sum_{i,j\geq 1}\left(ijp_{i+j}z\frac{\partial^2 
\Phi_h}{\partial p_i\partial p_j} +
                           ijp_{i+j}\frac{\partial\Phi_h }{\partial p_i}
                           \frac{\partial\Phi_h}{\partial p_j} +
          (i+j)p_ip_j\frac{\partial\Phi_h}{\partial p_{i+j}} \right).
\end{equation}}
\\

{\bf Remark}\quad  
By a combinatorial method,
I. P. Goulden and D. M. Jackson  have proven the above equation for the case $h=0.$ and 
gave  explicit formula of $\mu_{h,m}^{g,k}(\alpha)$ for $h=0$ and $ g=0,1,2$ (see $([GJ1,2,3], [GJV])$: 
\begin{eqnarray}
\mu_{0,m}^{0,k}(\alpha_1,\cdots,\alpha_m) &=&
	\frac{|c_\alpha|}{k!}(k+m-2)!k^{m-3}\prod_{i=1}^m
	\frac{\alpha_i^{\alpha_i}}{(\alpha_i-1)!}\\
\mu_{0,m}^{1,k}(\alpha_1,\cdots,\alpha_m) &=&
	\frac{|c_\alpha|}{24k!}(k+m)!\prod_{i=1}^m
        \frac{\alpha_i^{\alpha_i}}{(\alpha_i-1)!}(k^m-k^{m-1}-
        \sum_{i=2}^m(i-2)!e_ik^{m-i}),
\end{eqnarray}
where $c_\alpha$ is the conjugacy class of the symmetric group ${\cal 
S}_k$
on $k$ symbols indexed by the partition type $\alpha$ of $k$, and $e_i$ is 
the
$i$-th elementary symmetric function in $\alpha_1,\cdots,\alpha_m$.\\

For $h=0$, using the recursive formula (1.4) with initial value
$$\mu^{g,1}_{0,1}(1)=\left\{
   \begin{array}{lll}
    1&if&g=0\\
    0&if&g\geq 1\\
    \end{array}
    \right. ,$$  \\
we calculate
some  values of $\mu_{0,m}^{g,k}(\alpha)$ with the aid of {\bf Maple} for
$g=0,1,\cdots,5$, $k=3,4,5$, $m=1,2,\cdots,5$:\\
 \begin{center}
\begin{tabular}{|c|l|l|l|l|l|l|}\hline
$\alpha$&$\mu_{0,m}^{0,k}(\alpha)$&$\mu_{0,m}^{1,k}(\alpha)$&
$\mu_{0,m}^{2,k}(\alpha)$
&$\mu_{0,m}^{3,k}(\alpha)$&$\mu_{0,m}^{4,k}(\alpha)$&
$\mu_{0,m}^{5,k}(\alpha)$\\ \hline
(3)&1&9&81&729&6561&59049\\ \hline
(1,2)&4&40&364&3280&29524&265720\\ \hline
(1,1,1)&4&40&364&3280&29524&265720\\ \hline
(4)&4&160&5824&209920&7558144&272097280\\ \hline
(1,3)&27&1215&45927&1673055&60407127&2176250895\\ \hline
(2,2)&12&480&17472&629760&22674432&816291840\\ \hline
(1,1,2)&120&5460&206640&7528620&271831560&9793126980\\ \hline
(1,1,1,1)&120&5460&206640&7528620&271831560&9793126980\\ \hline
(5)&25&3125&328125&33203125&3330078125&333251953125\\ \hline
(1,4)&256&35840&3956736&409108480&41394569216&4156871147520\\ \hline
(2,3)&216&26460&2748816&277118820&27762350616&2777408868780\\ \hline
(1,1,3)&1620&234360&26184060&2719617120&275661886500&27700994510280\\ \hline
(1,2,2)&1440&188160&20160000&2059960320&207505858560&20803767828480\\ \hline
(1,1,1,2)&8400&1189440&131670000&13626893280&1379375197200&138543794363520\\ \hline
(1,1,1,1,1)&8400&1189440&131670000&13626893280&1379375197200&138543794363520\\ \hline
\end{tabular}
\end{center}
which coincide to the formulas (1.7) and (1.8) for $g$ = 0, 1, respectively.

To do similar calculation for $h>0$, we have to calculate some special Hurwitz number
$\mu_{h,m}^{g,k}(\alpha)$ for the case when $k+m-2+2g-2kh=0$,which we will discuss in another paper.
\\

{\bf Acknowledgements}  The first author would like to thank Prof.Kefeng Liu,
Prof. Zhengbo Qin and Prof. Yongbin Ruan for inviting him to visit their Universities
and for valuable discussions.The third author would like to thank Prof. Yongbin Ruan
and Prof. Renhong Wang for their honest helps.

\sect{2\quad  Relative GW-invariant}
\setcounter{section}{2}
\setcounter{equation}{0}
\\

Let $(M,\omega)$ be a real $2n$-dimension compact symplectic manifold 
with
symplectic form $\omega,$ and $Z^0,\ldots, Z^p$  symplectic submanifolds
of $M$
with real codimension $2$. Denote $Z=(Z^0,\cdots,Z^p)$. Let $\Sigma^g$
be a compact connected Riemman surface of genus $g\geq 0.$ Suppose
$A\in H_2(M,Z), k^i=\{k_1^i,\cdots,k_{l_i}^i\}$  a set of positive 
integers,
$i=0,\ldots,p$, denoted by $K=\{k^0,\ldots,k^p\}$. Consider moduli 
space ${\cal M} =
{\cal M}_{A,l}^{M,Z}(g,K)$ of pseudo-holomorphic maps $f:\Sigma^g\to M$ 
with marked points
$x_1,\cdots,x_l;$ $y_1^0,\cdots,y_{l_0}^0$; $\cdots$; 
$y_1^p,\cdots,y_{l_p}^p\in \Sigma^g$
such that $[f(\Sigma^g)]=A$, and $f$ is tangent to $Z^i$ at 
$y_1^i,\cdots,y_{l_i}^i$
with order $k_1^i,\cdots,k_{l_i}^i$, $i=0,\ldots,p$. Denote  
$x=(x_1,\cdots,x_l),$ $ y^i=(y_1^i,\cdots,y_{l_i}^i)$,
$y=(y^0,\cdots,y^p)$. Note that the intersection numbers $\#(A\cdot Z^i)$ 
are topological
invariants, and $\sum_{j=1}^{l_i}k_j^i=\#(A\cdot Z^i)$. Moreover, 
since $Z^i$ is a symplectic
submanifold, if $A$ can be expressed by the image of an nontrivial 
pseudo 
holomorphic map $f:\Sigma^g\to M$, the intersection number
$\#(Z_i \cdot A)\geq 0$.
Similarly to the Gromov-Uhlenbeck compactification for the  pseudo-holomorphic
maps, we comactify  ${\cal M}$ by $\overline{\cal M}=\overline{\cal
M}_{A,l}^{M,Z}(g,K)$, the space of relative stable maps (for details see [LR]).
 We have two natural maps:
 \begin{equation}
 \begin{array}{ll}
        \Xi_{g,l}:&\overline{\cal M}\to M^l\\
& (f,\Sigma^g, x, y, K)\longmapsto(f(x_1), \cdots, f(x_l))\\
\end{array}
\end{equation}
and
\begin{equation}\begin{array}{ll}
        p:& \overline{\cal M}\to Z_0\times\cdots\times Z_p\\
 & (f, \Sigma^g, x, y, K) \longmapsto ((f(y^0_1),
 \cdots, f(y^0_{l_0})),\cdots,
 (f(y^p_1),\cdots, f(y^p_{l_p}))).\\
  \end{array}
 \end{equation} \\

Roughly, the  relative GW-invariants are defined as
$$
        \psi_{A,g,l}^{M,Z}(\delta |\beta, K)= 
        \int_{\overline{\cal 
M}}\Xi_{g,l}^*\mathop{\Pi}\limits_{i}\delta_i\wedge p^*
        \mathop{\Pi}\limits_j\beta^j
$$
where $\delta=(\delta_1,\cdots,\delta_l)$, $\delta_i\in H^*(M,R)$,
$\beta=(\beta^0,\cdots,\beta^p)$, 
$\beta^j=(\beta_1^j,\cdots,\beta_{l_j}^j)$,
$\beta_i^j\in H^*(Z^j,\bf{R})$ for any $i$. For precise definition, see [LR].\\

If $\Sigma^g$ is not connected, suppose there exists $c$ connected 
components
$\Sigma^{g_1},\ldots,\Sigma^{g_c}$, then the genus $g$ is defined to be it's 
algebraic
genus, i.e., $g=\sum_{i=1}^c g_i-c+1$. Let ${\cal P}_x$ be the set 
of
all ordered partitions of $\{x_1,\ldots,x_l\}$ into $c$ parts.
Each $\pi=(\pi_1,\cdots,\pi_c)\in {\cal P}_x$
records which marked points $x=(x_1,\cdots,x_l)$ go on each 
components $\Sigma^{g_1}, \cdots,\Sigma^{g_c}.$ 
Similarly, we can define $\sigma^i=(\sigma_1^i,\cdots,\sigma_c^i)\in 
{\cal P}_{y^i}.$  Corresponding to the partition of $ y $, we define the
partition of $K$, i.e,  $\sigma^i=(\sigma_1^i,\cdots,\sigma_c^i)\in 
{\cal P}_{K^i}$ , and write
$\sigma_i=(\sigma_i^0, \cdots, \sigma_i^p)$, 
$\sigma=(\sigma_1,\cdots,\sigma_c)$. $\pi$, $\sigma$ induces a partition
of $\delta$, $\beta$, respectively. Denote the parameters over the component
$\Sigma^{g_i}$ by $x_{\pi_i}$, $y_{\sigma_i}$, $\delta_{\pi_i}$,
$\beta_{\sigma_i}$.
Suppose that $f_i : \Sigma \to M$ is a relative stable pseudo holomorphic map
such that $[f_i(\Sigma^{g_i} )]=A_i,$ and $\sum_{i=1}^{c}A_i=A.$  Then for
given $(A_1, \cdots, A_c),$ $\pi$ and $\sigma,$ the relative GW-invariant
$\psi_{A, g, l}^{M, Z, c}(\delta |\beta; K)$ is defined by
$$
 \psi_{A, g, l}^{M, Z, c}(\delta |\beta; K)(\pi, \sigma)=\prod_{i=1}^{c}
\psi_{A_i,g_i, l_i}^{M, Z_{\sigma_i}}(\delta_{\pi_i}|
\beta_{\sigma_i};K_{\sigma_i}).
$$
 \\

Consider the linearization of $\bar{\partial}$-operator
$$
 D_f=D_{\bar{\partial}_J(f)}:C^{\infty}(\Sigma, f^*TM)\to 
\Omega^{0,1}(f^*TM).
$$
If we choose a proper weighted Sobolev norm over 
$C^{\infty}(\Sigma, f^*TM)$
and $\Omega^{0,1}(f^*TM)$, we have the following, see [LR]\\

{\bf Lemma 2.1}\quad
{\it
$D_f$ is a Fredholm operator with index
\begin{equation}
{\rm Ind}\ (D_f) = 2C_1(M)A+(2n-6)(1-g)+2\sum_{i=0}^pl_i-
2\sum_{i=0}^p\sum_{j=1}^{l_i}k_j^i+2l
\end{equation}
and the  relative GW-invariant 
$\psi_{g,l}^{M,Z}(\delta|\beta;K)$
is defined to be zero unless
\begin{equation}
\sum_{i=1}^l \deg \delta_i + \sum_{i=0}^p\sum_{j=1}^{l_i} \deg 
\beta_j^i =
{\rm Ind}\ (D_f).
\end{equation}
}

Suppose that $H:M\to R$ is a proper periodic Hamiltonian function such 
that 
the Hamiltonian vector field $X^H$ generates a circle action. By adding 
a 
constant, we can assume that zero is regular value. Then $H^{-1}(0)$ is
a smooth submanifold preserved by circle action. The quotient 
$B=H^{-1}(0)/S^1$
is the famous symplectic reduction. Namely, $B$ has an induced 
symplectic
structure, so we can regard $B$ as a symplectic submanifold of $M$ with 
real
codimension 2. We cut $M$ along $H^{-1}(0)$. Suppose that we obtain two
disjoint components $M^{\pm}$ which have boundary $H^{-1}(0)$. We can 
collapse
the $S^1$-action on $H^{-1}(0)$ to obtain two closed symplectic 
manifolds
$\overline{M}^{\pm}$. This procedure is called as symplectic cutting, see 
[L], [LR].
Without  loss of generality, suppose $\overline{M}^+$ contains
$Z^+=(Z^0,\cdots,Z^q)$ as submanifolds and $\overline{M}^-$ contains
$Z^-=(Z^{q+1},\cdots,Z^p)$, $q\leq p$. There is a map
$$
        \pi:M\to \overline{M}^+\bigcup{}_B \overline{M}^-.
$$
It induces a homomorphism
$$
\pi^*: H^*(\overline{M}^+\cup_B \overline{M}^-,R)\to H^*(M,R).
$$
It was shown by Lerman [L] that the restriction of the symplectic 
structure
$\omega$ on $M^\pm$ such that $\omega^+|_B=\omega^-|_B$ is the induced 
symplectic form from symplectic reduction. By the Mayer-Vietoris 
sequence,
a pair of cohomology classes $(\delta^+,\delta^-)\in H^*(\overline{M}^+,R)
\otimes H^*(\overline{M}^-,R)$ with $\delta^+|_B=\delta^-|_B$
defines a cohomology class of $\overline{M}^+\cup_B\overline{M}^-$,
denoted by $\delta^+\cup_B\delta^-$.

Consider the moduli space
${\cal M}^+ ={{\cal 
M}^+}_{A^+,l^+}^{\overline{M}^+,Z^+,B,c^+}(g^+,K^+,\alpha^+)$
which consists of tuple $(\Sigma^{g^+}, x^+,$ $ y^+, e^+, K^+,
\alpha^+, f^+)$ with properties:
\begin{itemize}
\item  $\Sigma^{g^+}$ has $c^+$ connected components;
\item   $f^+:\Sigma^{g^+}\to \overline{M}^+$   is a  pseudo holomorphic map;
\item  $[f^+(\Sigma^{g^+})]=A^+$;
\item  $f^+$ is tangent to $Z^+=(Z_0,\cdots,Z_q)$ at
$y^+=(y^0,\cdots,y^q)$ with order $K^+=(k^0,\cdots,k^q)$; 
\item  $f^+$ is 
tangent to $B$ at $e^+=(e_1^+,\cdots,e_v^+)$ with order 
$\alpha^+=(\alpha_1^+,\cdots,\alpha_v^+)$.
\end{itemize}
Similarly, we can define
${\cal M}^- ={{\cal 
M}^-}_{A^-,l^-}^{\overline{M}^-,Z^-,B,c^-}(g^-,K^-,\alpha^-)$ which consists of
tuple $(\Sigma^{g^-}, x^-, y^-,$ $ e^-, K^-, \alpha^-, f^-)$.
 According to [LR], we can glue $f^+$ and $f^-$ to obtain 
a 
pseudo holomorphic map $f:\Sigma^g\to M$. A little more precisely, we glue 
$\overline{M}^+$
and $\overline{M}^-$ as above.
If $f^+$ and $f^-$ have same periodic orbits at each end , i.e, they have
same orders as they tangent to symplectic
submanifold
$B$, we can glue the maps $f^+$ and $f^-$ as $f^+\#f^-$ after 
gluing the
domain of Riemann surface $\Sigma^{g^+}$ and $\Sigma^{g^-}$, which is the 
connected sum of
$\Sigma^{g^+}$ and $\Sigma^{g^-}$.
 Then pertubating map $f^+\#f^-$,
we can get an
unique pseudo holomorphic map $f:\Sigma^{g}\to M$.
In our paper, we always require that  $\Sigma^{g^+}\#\Sigma^{g^-}$
is a connected Riemann surface. The
following  index
addition foumula is useful to our paper, \\

{\bf Lemma 2.2[LR]}\quad
$$
    {\rm Ind}(D_{f^+}) + {\rm Ind}(D_{f^-})=(2n-2)v+{\rm Ind}D_f. 
$$
\\

We also need a well known fact about genus of connected sum of Riemann surfaces:
\\

{\bf Lemma 2.3}\quad
{\it The following equality is satisfied:
\begin{equation}
        g=g^+ + g^- + v - 1
\end{equation}
where $g$ is the genus of $\Sigma^g$, $g^{\pm}$ is the algebraic genus of 
$\Sigma^{g^\pm}$,
$v$ is the number of end, i.e., the number of the points where we glue 
$\Sigma^{g^+}, \Sigma^{g^-}$. } \\

 According to 
theorem 5.8 of [LR], the   relative GW-invariant $\psi_{A,l}^{M,Z}(\delta | 
\beta;K)$ can be
expressed by the relative GW-invariants over each connected component.
 Precisely ,using the notations of [LR], suppose that
 ${\cal{C}}_{g,l,K}^{J,A}$ is the set of indices: 

(1) The combinatorial type of $ (\Sigma^{\pm}, f^{\pm}):\{A_i^{\pm},g_i^{\pm},
l_i^{\pm},K_i^{\pm},(\alpha^{\pm}_1,\cdots,\alpha^{\pm}_v)\},i=1,\cdots,v$,
$\sum\limits_{i=1}^{v}\alpha_i^\pm=\#(A\cdot B)$;

(2) A map $ \rho : \{e_1^+,\cdots, e_v^+\} \to \{e_1^-,\cdots,e_v^-\}$ ,where
$(e_1^{\pm},\cdots, e_v^{\pm})$ denote the puncture points of $\Sigma^{\pm}$,
satisfying

$(i)$ The map $\rho$ is one-to-one;

$(ii)$ If we identify $e_i^+$ and $\rho(e_i^+)$ , then
$\Sigma^+\bigcup \Sigma^-$ forms a connected closed Riemann surface of genus
$g$;

$(iii)$ $f^+(e_i^+)=f^-(\rho(e_i^+))$ and they have  same order of tangency;

$(iv)$ $((\Sigma^+,f^+),(\Sigma^-,f^-),\rho)$ represents the homology class
$A$.\\
For given $C\in {\cal{C}}_{g,l,K}^{J,A}$ suppose that $\pi^{\pm}_C$,
$\sigma^{\pm}_C$
are partitions of $x^{\pm}$, $y^\pm$, $e^\pm$, $\delta^\pm$, $\beta^\pm$,
$\alpha^\pm$
induced by $C$. Then we have the following:([LR] Lemma 5.4 and Theorem 5.8)\\

{\bf Lemma 2.4}\quad
{\it
${\cal{C}}_{g,l,K}^{J,A}$ is a finite set, and
\begin{eqnarray}
\psi_{A,g}^{M,Z}(\delta|\beta;K)=
\sum_{C \in{\cal{C}}_{g,l,K}^{J,A}}\psi_C(\delta|\beta;K),
\end{eqnarray}
where
\begin{eqnarray}
& &\psi_C(\delta|\beta;K)= \\
 & &\|\alpha\|\sum\delta^{I,J}
               \psi_{A^+,g^+,l^+}^{\overline{M}^+,Z^+,B,c^+}(\delta^+|
              \beta^+;\rho_{\mbox{}_I};K^+,\alpha)(\pi^+_C, \sigma^+_C)\cdot
{\psi}_{A^-,g^-,l^-}^{\overline{M}^-,Z^-,B,c^-}(\delta^-|
       \beta^-;\rho_{\mbox{}_J};K^-,\alpha)(\pi^-_C, \sigma^-_C),\nonumber\\  \nonumber
\end{eqnarray}
where $\|\alpha\|=\alpha_1\cdots\alpha_v$;
$\delta^{IJ}=\delta^{I_1J_1}\cdots
\delta^{I_vJ_v}$,  $\delta^{I_iJ_i}$ being the Kronecker symbol; and
$ \{\rho_1,\cdots,\rho_s\}$ is an orthonormal  basis of $H^*(B,\bf{R}), $
 $\rho_I=\{\rho_{I_1},\cdots,\rho_{I_v}\}$
$ \subset \{\rho_1,\cdots,\rho_s\}$,
$\rho_J=\{\rho_{J_1},\cdots,\rho_{J_v}\} \subset \{\rho_1,\cdots,\rho_s\}.$
 }
  \\

 For convenience in application, we will rewrite Lemma 2.4 in following steps:\\

 {\bf Step 1.} We divide $A$ into $A^+$ and $A^-$ such $ A = A^+\cup_B A^-.$\\

 {\bf Step 2.} Suppose $\Sigma^{g^\pm} $ have $a_i \geq 0 $ end points
 with order $i \in \{ 1, \cdots, \#(A \cdot B) \}$ such that $\sum\limits_{i}
 i\cdot a_i = \#(A\cdot B)$, and $g = g^+ + g^- + \sum\limits_{i}a_i - 1.$
 Denote $ a= (a_1, a_2, \cdots, ).$\\

{\bf Step 3.} Suppose that $ \tau^\pm = (\pi^\pm, \sigma^\pm) \in {\cal{ P}}_ {x^\pm}
\times {\cal{ P}}_{ y^\pm,  e^\pm}$ record which marked points in
$\{ x^\pm,$ $ y^\pm,$ $ e^\pm\} $ go on each component $\Sigma^{g^\pm_1},$
 $\cdots,$ $\Sigma^{g^\pm_{c^\pm}}$, satisfying:

 (1).$ g^{\pm}=\sum_{i=1}^{c^{\pm}} g_i^{\pm}-c^{\pm}+1, g_i^{\pm}\geq 0,
                i=1,\cdots,c^{\pm} $

 (2). $ f_i$ : $\Sigma^{g^+_i}\longrightarrow \overline{M}^\pm$ are relative stable
 holomorphic maps, and  $[f^+_i(\Sigma^{g^\pm_i})]=A^\pm_i,$
 $i=1, \cdots, c^\pm$  with
 $\sum\limits_{i=1}^{c^\pm}A_i^\pm=A^\pm.$\\

 Denote $\tau=(\tau^+, \tau^-).$ Note that $\tau^\pm$ induce a partition of
 $\delta^\pm$, $\beta^\pm$ , $a$ and $Z^\pm$.\\

  {\bf Step 4.} For given  $a$ and $\tau$,  we glue $\Sigma^{g^+}$ and
  $\Sigma^{g^-}$ in above manner such that $\Sigma^{g^+}\# \Sigma^{g^-}$ is a
  connected Riemann surface of genus $g$. However, for given such $a$ and
  $\tau$, we can glue
 $\Sigma^{g^+}$ and $\Sigma^{g^-}$ in many different ways
 such that $\Sigma^{g^+}\# \Sigma^{g^-}$ is a  connected Riemann surface of
 genus $g$. Denote  the number of different ways by $\kappa(a,\tau).$

Then we have the following
gluing formula for
the  relative GW-invariants :
\\

{\bf {Lemma 2.4$^{\prime}$}}\quad
{\it
\begin{eqnarray}
& &\psi_{A,g}^{M,Z}(\delta|\beta;K) =
 \sum\|a\| \cdot \delta^{IJ}\cdot\sum_{\tau}\kappa(a, \tau) \cdot  \\
   & &            \psi_{A^+,g^+,l^+}^{\overline{M}^+,Z^+,B,c^+}(\delta^+|
                \beta^+;\rho_{\mbox{}_I};K^+,a)(\pi^+, \sigma^+)\cdot\
{\psi}_{A^-,g^-,l^-}^{\overline{M}^-,Z^-,B,c^-}(\delta^-|
      \beta^-;\rho_{\mbox{}_J};K^-,a)(\pi^-, \sigma^-),\nonumber\\  \nonumber
\end{eqnarray}
where $\|a\|=1^{a_1}\cdot 2^{a_2}\cdot \cdots$,  and the first $\sum$ denotes that
we sum all possibility for
\\

\begin{equation}
\begin{array}{l}
        A = A^+ \bigcup_B A^-, \\
            g=g^++g^-+v-1, \\
            \#(A \cdot B)=\sum_i i \cdot a_i , \\
        \rho_I,\rho_J\subset \{\rho_1,\cdots,\rho_s\}. 
\end{array}
\end{equation}

}

\sect{3\quad  Relative GW-invariant over $\Sigma^h$}
\setcounter{section}{3}
\setcounter{equation}{0}
\\

In our case, $M$ is the Riemann surface $\Sigma^h$ with real two dimension, 
thus
$Z$ consists of points, which are the divisors of $M$. Since
$H_2(\Sigma^h,{\bf Z}) \cong{\bf Z}$,
denoted the generator by $H$, then the first Chern class $C_1(\Sigma^h)=(2-2h)H$.
Let $A=kH$. When we  say
$f\in {\cal M}_{A,l}^{\Sigma^h,Z}(g,K)$, we mean that $f:\Sigma^g\to \Sigma^h$
is a pseudo holomorphic map such that $[f(\Sigma^g)]=kH$ and there exists
marked points
$x,y\in \Sigma^g$, $f$ is tangent to $Z$ at $y$ with order $K$. Note that
$\sum_{j=1}^{l_i}k_j^i = \#(A\cdot Z_i)=\deg (f)=k$, then $k^i$ 
is a partition of
$k, i=0,\cdots,p$.
Moreover the  relative GW-invariant
$\psi_{A,g,0}^{\Sigma^h,Z}(|\beta;K)=0$ unless
\begin{equation}
\sum_{i=0}^p\sum_{j=1}^{l_i}\deg\beta_j^i = 2C_1(\Sigma^h)A+4(g-1)+
2\sum_{i=0}^{p} l_i-2\sum_{i=0}^p\sum_{j=1}^{l_i}k_j^i
\end{equation}
where $\beta_j^i\in H^*(Z^i,R)$, $j=1,\cdots,l_i$, $i=0,\cdots,p$. 
However,
since $Z^i$ is a point, $\deg \beta_j^i=0$. Thus we have 
\begin{equation}
 (2-2h)k+2(g-1)+
\sum_{i=0}^{p} l_i-\sum_{i=0}^p\sum_{j=1}^{l_i}k_j^i=0.
\end{equation}\\

{\bf Remark 3.1}\quad
{\it
The equality (3.2) is exactly the Riemann-Hurwitz formula. 
} \\

Suppose $\Sigma^g$ is connected, choose $l=0,
K=(2,1,\cdots,1;\cdots;2,1,\cdots,1;\alpha_1,\cdots,\alpha_m)$.
Then by  definition $\mu_{h,m}^{g,k}(\alpha)$ is just the  relative 
GW-invariant
$\psi_{A,g,0}^{\Sigma^h,Z}(|\beta;K)$.                                 
 
From (3.1), we derive $p=k+m-2kh-2+2g$, i.e., we have $k+m-2kh-2+2g$
double branch points over $\Sigma^h$, otherwise, $\mu_{h,m}^{g,k}(\alpha)=0$.

Now, we can prove  theorem A by symplectic cutting and the gluing formula
(2.8).
We perform the symplectic cutting over $\Sigma^h$ at $\infty$ in a small 
neighborhood
such that there is only one other double branch point $G$ in this 
neighborhood.
We have 
$$
\overline{M}^+=S^2,\quad \overline{M}^-=\Sigma^h
$$
It's easy to observe that $A^+=kH^{\prime}, A^-=kH$, where $H^{\prime}$ is the generator of
$H_2(S^2,{\bf Z}) \cong{\bf Z}$. We may consider 
dimension
condition equations:
\begin{equation}
\left\{
	\begin{array}{l}
        (k-m)+(2-1)+(k-v) =   2k-2+2g^+ \\
        g=g^+ + g^-+v-1
	\end{array}
\right.
\end{equation}

We first consider $\overline{M}^+=S^2$. The map
$f^+:\Sigma^{g^+}\to \overline{M}^+$ branches
at only three points: infinity, the fixed double branch point $G$ and
the symplectic reduction point $B$.
Suppose $\Sigma^{g^+}=\cup_{i=1}^{c^+}\Sigma^{g_i^+}$, i.e., 
$\Sigma^{g^+}$
has $c^+$-connected components.
Suppose the holomorphic map $u_i^+:\Sigma^{g_i^+}\to\overline{M^+}$ has degree 
$k_i^+$. It is obvious that $k_i^+\leq k$.
If $\Sigma^{g_i^+}$ contains a double ramification point,
we have from Riemann-Hurwitz formula
over this component that
\begin{equation}
k_i^+-v_i^++2-1+k_i^+-m_i^+=2k_i^+-2+2g_i^+
\end{equation}
where $v_i^+\geq 1, m_i^+\geq 1$ is the number of ramification at the
symplectic reduction point $B$ and infinity, respectively. Note that
the geometric genus $0\leq g_i^+\leq g$, so we have two cases:
$(v_i^+,m_i^+,g_i^+)=(1,2,0)$, or $(2,1,0)$. By the same reason, if 
the 
component $\Sigma^{g_i^+}$ doesn't contain any double ramification
point, we have
one case $(v_i^+,m_i^+,g_i^+)=(1,1,0)$. Note that 
$\sum_{i=1}^{c^+}m_i^+=m$,
we have $v=\sum_{i=1}^{c^+}v_i^+=m-1$, or $m+1$, correspondingly,
$c^+=m-1$, or, $m$. In sum, we have proven  \\

{\bf Lemma 3.2}\quad
{\it
For $\overline{M}^+=S^2$, the holomorphic map $f_i:\Sigma^{g_i^+}\to 
\overline{M}^+$ has
one of the following branch types

(1)\quad $(\alpha_i;1,1,\cdots,1;\alpha_i), 1\leq i\leq m$

(2)\quad $(\alpha_k,\alpha_l; 2,1,\cdots,1;\alpha_k+\alpha_l), 1\leq 
k<l\leq m$

(3)\quad $(\alpha_i;2,1,\cdots,1;\rho,\alpha_i-\rho), 1\leq i\leq m, 
1\leq\rho
	  \leq \left[\frac{\alpha_i}{2}\right]$\\
at infinity , the fixed branch point $G$ and the symplectic reduction
point $B$ respectively.
 }
\\

If $v=m-1$, then $c^+=m-1, g^+=\sum_{i=1}^{c^+}
g_i^+-c^++1=2-m$.
 Substituting into $g=g^++g^-+v-1$, we have $g^-=g$. Note that
$g^-=\sum_{i=1}^{c^-}g_i^{-}-c^-+1$,
$0\leq \sum g_i^-\leq g$, $g_i^-\geq 0, c^-\geq 1$, we have only one 
case: $c^-=1, g^-=g$.
If  $v=m+1,$ then $ c^+=m$.
By the same reason, we have two cases:
 $c^-=1, g^-=g-1$ and
 $c^-=2,$ $g^-= g_1^-+g_2^--2+1=g-1,$ $ g_1^+\geq 0, g_2^-\geq 0$. In sum, we have
 \\

{\bf Lemma 3.3}\quad
{\it
The genus $g^-$ and the number $c^-$ of connected components of
$\Sigma^-$ are one of the following cases:
 
(i) $c^-=1, g^-=g$;

(ii)  $c^-=1, g^-=g-1$;

(iii)  $c^-=2$,  $g^-= g_1^-+g_2^--1=g-1$, $ g_1^-\geq 0,  g_2^-\geq 0$. }
\\

  Regarding the symplectic reduction point $B\in \overline{M}^-$ as infinity, we get many
new almost simple ramified covering maps $f_i^-:\Sigma^{g_i^-}\to 
\overline{M}^-$. However in
any above case, the holomorphic map $f_i^-:\Sigma^{g_i^-}\to \overline{M}^-$ 
has either
strict smaller number of ramification points  at infinity, or strict smaller 
degree, or
strict smaller genus than the holomorphic map $f:\Sigma^g\to M=\Sigma^h$. 
Thus if we
have known the  relative GW-invariant in $\overline{M}^+$, by 
Lemma 2.4, we can get a recursive formula for $\mu_{h,m}^{g,k}(\alpha)$.
We also need the following lemma.
 \\

{\bf Lemma 3.4}\quad
{\it
Let 
$\theta=(\alpha_1,\cdots,\hat{\alpha_i},\cdots,\hat{\alpha_j},\cdots,
\alpha_m,\alpha_i+\alpha_j)\in J(\alpha)$,
$\omega=(\alpha_1,\cdots,\hat{\alpha_i}$,
$\cdots,\alpha_m,\rho$, $\alpha_i-\rho)
\in C(\alpha)$ as in the introduction, then for $\overline{M}^+$,
the product $\psi_J(\alpha,\theta)$ of the relative
GW-invariants of $(m-1)$ connected components is
\begin{equation}
\psi_J(\alpha,\theta) 
        = \left\{
		\begin{array}{ll}
                \frac{1}{\alpha_1}\cdots\hat{\frac{1}{\alpha_i}}\cdots
                \hat{\frac{1}{\alpha_j}}\cdots
		\frac{1}{\alpha_m}  & \mbox{ if } \alpha_i\not=\alpha_j\\
                \frac{1}{\alpha_1}\cdots\hat{\frac{1}{\alpha_i}}\cdots
                \hat{\frac{1}{\alpha_j}}\cdots
		\frac{1}{\alpha_m}\cdot\frac{1}{2}
                  & \mbox{ if } \alpha_i=\alpha_j
		\end{array}
        \right. \nonumber
\end{equation}
and the product $\psi_C(\alpha, \omega)$ of the relative
GW-invariants of $m$ connected components  is
\begin{equation}
\psi_C(\alpha,\omega)        
        = \left\{
		\begin{array}{ll}
                \frac{1}{\alpha_1}\cdots\hat{\frac{1}{\alpha_i}}\cdots
                \frac{1}{\alpha_m}  & \mbox{ if } \rho\not=\alpha_i-\rho,\\
                \frac{1}{\alpha_1}\cdots\hat{\frac{1}{\alpha_i}}\cdots
		\frac{1}{\alpha_m}\cdot\frac{1}{2}
                  & \mbox{ if } \rho=\alpha_i-\rho.
		\end{array}
        \right. \nonumber
\end{equation}
} 

{\bf Proof}\quad
By the definition of the  relative GW-invariant and Lemma 
3.2,
we only need to calculate the connected  relative 
GW-invariant of
two type:
\begin{equation}
{\large
\begin{array}{l}
Q_1 =  \psi_{kH,0,0}^{S^2,pt,pt,pt}(|pt,pt,pt;k;1,1,\dots,1;k),
        k\geq 1\\[12pt]
Q_2 =  \psi_{kH,0,0}^{S^2,pt,pt,pt}
        (|pt,pt,pt;k;2,1,\cdots,1;\rho,k-\rho),\ \rho\geq 1,
\end{array}}
\end{equation}
where the ``pt'' in the bracket records the point homology Poincare dual 
to
the generator $E\in H^0(pt,R)$, and the others correspond to $Z$.

Regarding a holomorphic map $f:S^2\to S^2$ as a meromorphic function 
over 
Riemann plane, we  write $f\in$
$\overline{{\cal M}}_{kH,0}^{S^2, pt, pt, pt}(0; K; 1,1,\cdots, 1;k)$
 in the form $F_1:C\to C$,
$F_1(x)=\frac{\alpha_0(x-y^1)^k}{(x-y^2)^k},x\in C$, where 
$y_1\not=y_2\in C$
are $k$-ramification points. Without loss of generality, we choose zero and
infinity as k-ramification points, and send 1 to 1 ,thus there exists
 an uniqe
solution $F_1(x)=x^k$. However we have conformal
transformation $\pi_i:C \to C ,
\pi_i(x)=e^{\frac{2\pi i}{k}}x, i=0,\cdots,k-1, $ such that
$F_1(x)=F_1\circ\pi_i(x)$
 , i.e, there exists a finite group ${\bf Z}_k$ that acts on
$\overline{{\cal M}}_{kH,0}^{S^2, pt, pt, pt}(0; K; 1,1,\cdots, 1;k)$, 
 thus
 \begin{equation}
  Q_1=\frac{1}{k}.
 \end{equation}

By the same reason, we write $f\in$
$\overline{{\cal M}}_{kH,0}^{S^2, pt, pt, pt}(0; K; 2,1,
\cdots, 1;\rho, k-\rho)$
in form $F_2:C\to C$,
$$F_2(x)=\frac{\alpha_0(x-y_1^1)^\rho(x-y_2^1)^{k-\rho}}{(x-y^2)^k},$$
$\alpha_0 \neq 0 , x\in C$,
where $y_1^1\not\not=y_2^1$, $y^2\in C$ are $\rho, k-\rho$, 
$k$-ramification point, respectively. Suppose 1,2 and zero
are $\rho,k-\rho,k$-
ramification point, respectively. Then $F_2=\frac{\alpha_0(x-1)^\rho(x-2)   
 ^{k-\rho}}{x^k}$. Since $ F_2$ has a double ramification
 point $x$ at a given point
 , for instance at 1 . Then $x \neq 0,1,2$. We have following
 equation
\begin{equation}
      \left\{
      \begin{array}{lll}
      F_2(x)&=&1,\\
       F_2^{\prime}(x)&=&0.\\
      \end{array}
      \right.
\end{equation}
Solving(3.8), we have unique solution
$F_2=\frac{\alpha_0(x-1)^\rho(x-2)^{k-\rho}}
{x^k}$, where $\alpha_0=\frac{(2k)^k}{\rho^{\rho}
(2\rho-2k)^{k-\rho}}.$
 However, if $\rho=k-\rho$, we have conformal transformation
  $\pi:C \to C,\pi(x)=\frac
{2x}{3x-2}$,such that $F_2(x)=F_2 \circ\pi(x)$.
Since $\pi \circ \pi ={\bf 1}$,
there exists a finite group ${\bf Z}_2$ that acts on
$\overline{{\cal M}}_{kH,0}^{S^2, pt, pt, pt}(0; K; 2,1,\cdots, 1;
\rho, k-\rho)$, thus
\begin{equation}
Q_2 = \left\{
	\begin{array}{lll}
        \frac{1}{2} & \mbox{if} & \rho=k-\rho,\\
        1 & \mbox{if} & \rho\not=k-\rho.\\
        \end{array}
\right.
\end{equation}
We complete the Lemma 3.4. $\Box$
                            \\

Now, we prove theorem A:\\

{\bf  Theorem A}\quad
{\it 
Hurwitz number $\mu_{h,m}^{g,k}(\alpha)$ can be determined by a recursive
formula:
\begin{eqnarray}
\mu_{h,m}^{g,k}(\alpha)
     &=& \sum_{\theta\in 
J(\alpha)}\mu_{h,m-1}^{g,k}(\theta)\cdot I_1(\theta)
     + \sum_{\omega\in 
C(\alpha)}\mu_{h,m+1}^{g-1,k}(\omega)\cdot I_2(\omega)\\
     & & + \sum_{\omega\in C(\alpha)}
           \sum_{\mbox{\scriptsize $\begin{array}{c} g_1+g_2=g \\ 
g_1,g_2\geq 0
                                          \end{array}$}} 
           \sum_{\mbox{\scriptsize $\begin{array}{c} k_1+k_2=k \\ 
k_1,k_2\geq 1
					  \end{array}$}}
           \sum_{\mbox{\scriptsize $\begin{array}{c} m_1+m_2=m+1\\ 
m_1,m_2\geq 1
					  \end{array}$}}
           \sum_{\pi\in {\cal P}_{\omega}} \nonumber\\
     & &    {k+m-2kh-3+2g\choose k_1+m_1-2k_1h-2+2g_1}
            \mu_{h,m_1}^{g_1,k_1}(\omega_{\pi_1}) \cdot
            \mu_{h,m_2}^{g_2,k_2}(\omega_{\pi_2})\cdot
             I_3(\pi). \nonumber
\end{eqnarray}
where $m_{i}=l(\omega_{\pi_i}), i=1,2.$
}

{\bf Proof}\quad
For a positive integer $b$ and an ordered positive integer tuple
$\beta=(\lambda_1, \cdots, \lambda_t)$, we define an integer
$\varphi(\beta, b)=\#\{\lambda\in\beta|\lambda=b\}.$
According to  Lemma 2.4$^\prime$ and Lemma 3.3, we have
\begin{eqnarray}
  &\mu_{h,m}^{g,k}(\alpha)& \nonumber\\
     &=& \sum_{\theta\in J(\alpha)}\mu_{h,m-1}^{g,k}(\theta)\cdot
        \psi_J(\alpha,\theta)\cdot\|\theta\|\cdot\varphi(\theta,
        \alpha_i+\alpha_j)                 \nonumber\\
    &&  + \sum_{\omega\in C(\alpha)}\mu_{h,m+1}^{g-1,k}(\omega)\cdot
                \psi_C(\alpha,\omega)\cdot\|\omega\|
            \cdot\varphi(\omega, \rho)\cdot(\varphi(\omega, \alpha_i-\rho)
                -\delta^{\rho,\alpha_i-\rho})          \nonumber\\
     && + \sum_{\omega\in C(\alpha)}
           \sum_{\mbox{\scriptsize $\begin{array}{c} g_1+g_2=g\\ 
g_1,g_2\geq 0
					  \end{array}$}}
           \sum_{\mbox{\scriptsize $\begin{array}{c} k_1+k_2=k\\ 
k_1,k_2\geq 1
					  \end{array}$}}
           \sum_{\mbox{\scriptsize $\begin{array}{c} m_1+m_2=m+1\\ 
m_1,m_2\geq 1
					  \end{array}$}}
         \sum_{\pi\in {\cal P_{\omega}}} {k+m-2kh-3+2g\choose k_1+m_1-2k_1h-2+2g_1} \\           
 & &    \cdot \mu_{h,m_1}^{g_1,k_1}(\omega_{\pi_1})\cdot
      \varphi(\omega_{\pi_1}, \rho)
   \cdot\mu_{h,m_2}^{g_2,k_2}(\omega_{\pi_2})\cdot\varphi(\omega_{\pi_2},
   \alpha_i-\rho) \cdot\psi_C(\alpha,\omega)\cdot\|\omega\|  \nonumber
\end{eqnarray}
where $\|\theta\|=\alpha_1\cdots\hat{\alpha}_i\cdots\hat{\alpha}_j
\cdots\alpha_m(\alpha_i+\alpha_j);$         
 $\|\omega\|=\alpha_1\cdots\hat{\alpha}_i
\cdots\alpha_m\rho(\alpha_i-\rho);$ 
 $\delta^{\rho, \alpha_i-\rho}$ is the Kronecker symbol;
 $m_i=l(\omega_{\pi_i}), i=1,2$;
the factor ${k+m-2kh-3+2g\choose k_1+m_1-2k_1h-2+2g_1}$ comes from  the fact that we can
choose
$k_1+m_1-2k_1h-2+2g_1$  double ramification points over the component 
$\Sigma^{g_1}$
from $k+m-2kh-3+2g$ double ramification points.
Substituting (3.4), (3.5) into  (3.11), we get  (3.10).
 $\Box$
\\

{\bf Remark}\quad For $h=0$, the initial value of our recursive formula is
$$\mu^{g,1}_{0,1}(1)=\left\{
   \begin{array}{lll}
    1&if&g=0\\
    0&if&g\geq 1\\
    \end{array}
    \right. ,$$  

For $h>0$,the initial value  is not only
$$\mu^{g,1}_{h,1}(1)=\left\{
   \begin{array}{lll}
    1&if&g=h\\
    0&if&g\geq h+1\\
    \end{array}
    \right., $$\\
but also some special Hurwitz number $\mu^{g,k}_{h,m}(\alpha)$ for the case
when $k+m-2kh-2+2g=0$, which we will discuss in another paper.

\mbox{}\\[10pt]
\noindent{\Large\bf References}
\begin{description}
\leftmargin = 1.6cm
\item[\mbox{\makebox[1cm][l]{[A]}}]
  V.I.Arnol'd, Topological classification of trigonometric polynomial
  and combinatorics of graphs with an equal number of vertices and edges,
  Func.Ann. and its Appl.,30 No.1(1996),1--17.
\item[\mbox{\makebox[1cm][l]{[GJ1]}}]
  I.P. Goulden \& D.M. Jackson, A proof of a conjecture for the
  number of ramified covering of the sphere by the torus, preprint.
\item[\mbox{\makebox[1cm][l]{[GJ2]}}]
  ---, The number of ramified covering of the sphere by the double
  torus, and a general form for higher genera, preprint.
\item[\mbox{\makebox[1cm][l]{[GJ3]}}]
  ---, Transitive factorisations into transpositions and holomorphic
  mapping on the sphere, Proceeding of AMS Vol(125) No. 1(1997), 51--60.
\item[\mbox{\makebox[1cm][l]{[GJV]}}]
  I.P.Goulden,D.M.Jckson \& A.Vainshtain, The number of ramified covering
  of the sphere by the torus and surface of higher  genera, AG/9902125.
\item[\mbox{\makebox[1cm][l]{[CT]}}]
  M.Crescimanno \& W.Taylor ,Large $N$ phases of chiral $QCD_2$, Nuclear Phys. B,
  437 No.1(1995),3--24.
\item[\mbox{\makebox[1cm][l]{[D]}}]
  J.D\'{e}nes,The representation of a permutation as the product of a
  minimal number of transpositions and its connection with the  theory of graphs
  ,Publ.Math.Inst.Hungar.Acad.Soc., 4(1959),63--70.
\item[\mbox{\makebox[1cm][l]{[H]}}]
  A. Hurwitz. Ueber Riemann'sche Fl\"{a}chen mit gegebenen
  Verzweigungspunkten, Math. Ann, 39(1891), 1--60.
\item[\mbox{\makebox[1cm][l]{[IP]}}]
  E. Ionel \& T. Parker, Gromov-Witten invariants of symplectic sums.
  Prepint. Math.sg/9806013.
\item[\mbox{\makebox[1cm][l]{[L]}}]
  E. Lerman, Symplectic cuts, Math. Research Let. 2(1985) 247--258.
\item[\mbox{\makebox[1cm][l]{[LR]}}]
  An-Min Li \& Yongbin Ruan, Symplectic surgery and Gromov-Witten
 invariants of Calabi-Yau 3-folds I, preprint. Math.alg-geom/9803036.
\item[\mbox{\makebox[1cm][l]{[LLY]}}]
 B. Lian, K. Liu \& S.-T. Yau, Mirror principle I, Asian J. Math. 1(1997), 729-763.
\item[\mbox{\makebox[1cm][l]{[V]}}]
 R.Vakil, Recursions,formulas,and graph-theoretic interpretation of
 ramified coverings of the sphere by surface of genus 0 and 1, CO/9812105.
\end{description}
\end{document}